\documentstyle[11pt,twoside]{article}
\title{The asymptotic expansion of a generalisation of the Euler-Jacobi series}
\author{\sc R. B.\ Paris \\
{\em Division of Computing and Mathematics}, \\
{\em University of Abertay Dundee, Dundee DD1 1HG, UK}
}
\begin{document}
\def\f#1#2{\mbox{${\textstyle \frac{#1}{#2}}$}}
\def\dfrac#1#2{\displaystyle{\frac{#1}{#2}}}
\def\boldal{\mbox{\boldmath $\alpha$}}
{\newcommand{\Sgoth}{S\;\!\!\!\!\!/}
\newcommand{\bee}{\begin{equation}}
\newcommand{\ee}{\end{equation}}
\newcommand{\lam}{\lambda}
\newcommand{\ka}{\kappa}
\newcommand{\al}{\alpha}
\newcommand{\th}{\theta}
\newcommand{\fr}{\frac{1}{2}}
\newcommand{\fs}{\f{1}{2}}
\newcommand{\g}{\Gamma}
\newcommand{\br}{\biggr}
\newcommand{\bl}{\biggl}
\newcommand{\ra}{\rightarrow}
\newcommand{\mbint}{\frac{1}{2\pi i}\int_{c-\infty i}^{c+\infty i}}
\newcommand{\mbcint}{\frac{1}{2\pi i}\int_C}
\newcommand{\mboint}{\frac{1}{2\pi i}\int_{-\infty i}^{\infty i}}
\newcommand{\gtwid}{\raisebox{-.8ex}{\mbox{$\stackrel{\textstyle >}{\sim}$}}}
\newcommand{\ltwid}{\raisebox{-.8ex}{\mbox{$\stackrel{\textstyle <}{\sim}$}}}
\renewcommand{\topfraction}{0.9}
\renewcommand{\bottomfraction}{0.9}
\renewcommand{\textfraction}{0.05}
\newcommand{\mcol}{\multicolumn}
\date{}
\maketitle
\pagestyle{myheadings}
\markboth{\hfill \sc R. B.\ Paris  \hfill}
{\hfill \sc  Poisson-Jacobi-type transformation\hfill}
\begin{abstract}
We consider the asymptotic expansion of the sum
\[S_p(a;w)=\sum_{n=1}^\infty \frac{e^{-an^p}}{n^{w}}\]
as $a\rightarrow 0$ in $|\arg\,a|<\fs\pi$ for arbitrary finite $p>$ and $w>0$. 
Our attention is concentrated mainly on the case when $p$ and $w$ are both even integers, where the expansion consists of a {\it finite} algebraic expansion together with a sequence of increasingly subdominant exponential expansions. This exponentially small component produces a transformation for $S_p(a;w)$ analogous to the well-known Poisson-Jacobi 
transformation for the sum with $p=2$ and $w=0$.  Numerical results are given to illustrate the accuracy of the expansion obtained.

\vspace{0.4cm}

\vspace{0.3cm}

\end{abstract}

\vspace{0.3cm}

\noindent $\,$\hrulefill $\,$

\vspace{0.2cm}

\begin{center}
{\bf 1. \  Introduction}
\end{center}
\setcounter{section}{1}
\setcounter{equation}{0}
\renewcommand{\theequation}{\arabic{section}.\arabic{equation}}
We consider the asymptotic expansion of the sum
\bee\label{e11}
S_p(a;w)=\sum_{n=1}^\infty \frac{e^{-an^p}}{n^w}
\ee
as the parameter $a\ra 0$ in $|\arg\,a|<\fs\pi$, where $p>0$ and, for convenience, $w$ will be supposed throughout to be real
and positive.  When $w=0$, this sum is known as the Euler-Jacobi series and
when $a=0$ then $S_p(0;w)$ reduces to the Riemann zeta function $\zeta(w)$ (provided $\Re (w)>1$). Consequently, the series in (\ref{e11}) can also be viewed as a smoothed Dirichlet series for $\zeta(w)$.

The asymptotics of $S_p(a;w)$ as $a\ra 0+$ for $p$ a rational fraction and $w<0$ was considered by Ramanujan and is discussed in \cite[Chapter 15]{BR}. It was shown that the expansion in this case consisted of an asymptotic sum involving the Riemann zeta function.
A hypergeometric function approach for the Euler-Jacobi series when $p$ is a rational fraction and $w=0$ has been discussed at length in the monograph \cite{K}, where great emphasis was placed on obtaining exponentially small expansions. The case $w=0$ and arbitrary $p>0$, thus generalising the work in \cite{K}, 
has been investigated in \cite[\S 8.1]{PK} also for $a\ra 0+$ using both a Mellin-Barnes integral approach and also a saddle point analysis of a Laplace-type integral representation. A similar Mellin-Barnes integral approach for this latter case has also been independently considered in \cite{B}. The saddle point approach is instructive for understanding the appearance of exponentially small terms in the expansion as the parameter $p$ increases. It was established in \cite[\S 8.1]{PK} that an additional exponentially small contribution appears when $p=2+4k$, $k=0, 1, 2, \ldots\,$.  Indeed, the appearance 
of the first exponentially small expansion as $p$ passes through the `classical' value $p=2$ was demonstrated to be associated with a Stokes phenomenon \cite[\S 8.1.7]{PK}.  

In the case $p=2$, $w=0$ the behaviour of the sum $S_2(a;0)$ as $a\ra 0$ can be obtained from
the classical Poisson-Jacobi transformation given by
\bee\label{e12}
S_2(a;0)=\sum_{n=1}^\infty e^{-an^2}=\frac{1}{2}\sqrt{\frac{\pi}{a}}-\frac{1}{2}+\sqrt{\frac{\pi}{a}} \sum_{n=1}^\infty e^{-\pi^2n^2/a}
\ee
valid for all values of $a$ in $|\arg\,a|<\fs\pi$.
This well-known transformation relates a sum of Gaussian exponentials involving the parameter $a$ to a similar sum with parameter $\pi^2/a$. In the small-$a$ limit,
the convergence of the sum on the left-hand side becomes slow, whereas the sum on the right-hand side converges rapidly. Various proofs of (\ref{e12}) exist in the literature; see, for example, 
\cite[p.~120]{PK}, \cite[p.~60]{T1} and  \cite[p.~124]{WW}.

The dominant asymptotic expansion of $S_p(a;w)$ as $a\ra 0$ in $|\arg\,a|<\fs\pi$ for general $p>0$ and $w>0$ is relatively straightforward and is found to consist, in general, of a single term proportional to $a^{(w-1)/p}$, together with a series in ascending powers of $a$ with coefficients involving the Riemann zeta function (the algebraic expansion). When $0<p\leq 1$, the expansion is convergent and the result is exact; when $p>1$ the expansion is asymptotic as $a\ra 0$.
The most interesting case arises when $p$ and $w$ are both even integers. 
The above-mentioned algebraic expansion then terminates after a finite number of terms, and it becomes essential for accurate estimation to also include a subdominant sequence of exponentially small expansions. This exponentially small component produces a transformation
for $S_p(a;w)$ analogous to the Poisson-Jacobi transformation in (\ref{e12}), but valid as $a\ra 0$ in $|\arg\,a|<\fs\pi$. This similarly involves a finite sequence of series similar to (\ref{e12}) with $a$ in the exponential replaced by an inverse power of $a$, but with each term decorated by an asymptotic series in ascending powers of $a^{1/(p-1)}$. 

The approach we employ in this paper is based on a Mellin-Barnes integral representation for $S_p(a;w)$ and is similar to that described in \cite[\S 8.1.4]{PK}. An algorithm for the determination of the coefficients in the exponentially small asymptotic series is described.
The case $p=2$ when $w$ is an even integer has been recently discussed in \cite{P}, where the coefficients in the decorating asymptotic series can be given in closed form.
An application of the series when $p=2$, with $w=2$ and $w=4$, has arisen in the geological problem of thermochronometry  in spherical geometry \cite{WFK}.
\vspace{0.6cm}

\begin{center}
{\bf 2. \ An expansion for $S_p(a;w)$ as $a\ra 0$ when $w,\ p\neq 2,4, \ldots$}
\end{center}
\setcounter{section}{2}
\setcounter{equation}{0}
\renewcommand{\theequation}{\arabic{section}.\arabic{equation}}
We examine the expansion of the series $S_p(a;w)$ defined in (\ref{e11}) as $a\ra 0$ in the sector $|\arg\,a|<\fs\pi$, where $p>0$ and  for simplicity in presentation we shall assume throughout real values of $w>0$.
The case $p=1$, $w=0$ may be excluded from our consideration since the series in this case is summable as a geometric progression.

Our starting point is the well-known Cahen-Mellin integral (see, for example, \cite[\S 3.3.1]{PK})
\bee\label{e21}
z^{\alpha} e^{-z}=\frac{1}{2\pi i}\int_{c-\infty i}^{c+\infty i} \g(\alpha-s) z^{s}ds \qquad(z\neq 0,\ |\arg\,z|<\fs\pi),
\ee
where $c<\Re (\alpha)$ so that the integration path passes to the left of all the poles of $\g(\alpha-s)$ situated at $s=k+\alpha$ ($k=0, 1, 2, \ldots$).
Then, it follows that
\begin{eqnarray}
S_p(a;w)&=&\sum_{n=1}^\infty \frac{e^{-an^p}}{n^w}
=\sum_{n=1}^\infty\frac{n^{-w}}{2\pi i}\int_{-c-\infty i}^{-c+\infty i} \g(-s) (an^p)^{s}ds\nonumber\\
&=&\frac{1}{2\pi i}\int_{-c-\infty i}^{-c+\infty i} \g(-s) \zeta(w-ps)a^{s}ds\qquad (|\arg\,a|<\fs\pi),\label{e21a}
\end{eqnarray}
upon reversal of the order of summation and integration, which is justified when $c>\max\{0,(w-1)/p \}$, and evaluation of the inner sum in terms of the Riemann zeta function.

The integrand in (\ref{e21a}) possesses simple poles at $s=k$ ($k=0, 1, 2, \ldots$) and $s=s_0\equiv (w-1)/p$, except if $s_0=M$ (that is, $w=pM+1$), where $M$ is a non-negative integer, when the pole at $s=s_0$ is double. 
The residue at the double pole is obtained by making use of the fact that $\zeta(s)\simeq 1/(s-1)+\gamma$ in the neighbourhood of $s=1$, where $\gamma$ is Euler's constant, to find
\[\frac{(-a)^M}{M!} \bl\{\gamma-\frac{1}{p}\log\,a+\frac{1}{p}\,\psi(M+1)\br\}\qquad (M=0, 1, 2, \ldots),\]
where $\psi(x)$ is the logarithmic derivative of the gamma function.
The case when $w$ and $p$ are even positive integers requires a separate treatment which is discussed in Section 3.
\vspace{0.3cm}

\noindent 2.1\ {\it The case $0<p<1$}
\vspace{0.2cm}

\noindent
We first consider the case $0<p<1$. The integration path in (\ref{e22}) can be made to coincide with the imaginary $s$-axis together with a suitable indentation to lie to the left of the poles at $s=0$ and $s=s_0$ (when $0<w<1$). Then use of the functional relation for $\zeta(s)$ given by \cite[p.~269]{WW}
\bee\label{e23a}
\zeta(s)=2^s\pi^{s-1} \zeta(1-s) \g(1-s) \sin \fs\pi s
\ee
shows that the integrand can be written as
\[(2\pi)^w\zeta(1-w+ps)\,\frac{\g(1-w+ps)}{\g(1+s)}\,\frac{\sin \fs\pi(ps-w)}{\sin \pi s}\,\frac{a^s}{(2\pi)^{ps}}.\]
With $s=Re^{i\theta}$, where $R\ra\infty$ and is chosen so that the arc passes between the poles on the positive real axis, the logarithm of the dominant real part of the integrand is controlled by
\[(p-1)R\cos \theta\,\log\,R+O(R).\]
When $|\theta|<\fs\pi$, this last expression tends to $-\infty$ when $0<p<1$. Consequently, the integration path can be bent back to enclose the poles of the integrand to yield the convergent result
\bee\label{e22}
S_p(a;w)=J_p(a;w)+\sum_{k=0}^{\infty}{}^{\!'}\frac{(-)^k}{k!}\,\zeta(w-kp)a^k\qquad (0<p<1),
\ee
where 
\[J_p(a;w)=\left\{\begin{array}{ll}\dfrac{1}{p}\g\bl(\dfrac{1-w}{p}\br) a^{(w-1)/p} & (w\neq pM+1) \\
\\
\dfrac{(-a)^M}{M!} \{\gamma-\dfrac{1}{p}\log\,a+\dfrac{1}{p}\,\psi(M+1)\} & (w=pM+1) \end{array}\right.\]
and the prime on the sum over $k$ denotes the omission of the term corresponding to $k=M$ when $w=pM+1$.

When $p=1$, the sum $S_1(a;w)$ is given by (\ref{e22}) but now the sum over $k$ on the right-hand side converges when $|a|<2\pi$; see \cite[\S 4.2.2]{PK} for details. If we let $p=1$, $w=0$ then use of the facts that $\zeta(1-2k)=-B_{2k}/(2k)$ and $\zeta(-2k)=0$, where $B_{2k}$ are even-order Bernoulli numbers, shows that
\[S_1(a;0)=\frac{1}{a}-\frac{1}{2}+\frac{1}{a}\sum_{k=1}^\infty \frac{B_{2k}}{(2k)!}\,a^{2k}\qquad (|a|<2\pi),\]
which correctly reduces to the trivial summation $1/(e^a-1)$  by application of \cite[Eq.~(24.2.1)]{DLMF}.
\vspace{0.3cm}

\noindent 2.2\ {\it The case $p>1$}
\vspace{0.2cm}

\noindent
When $p>1$, the integration path in (\ref{e21a}) cannot be bent back over the poles and we proceed in a similar manner to that described for the case $w=0$ in \cite[\S 8.1.4]{PK}.
Consider the integral taken round the rectangular contour with vertices at $-c\pm iT$, $c'\pm iT$, where $c'>0$. The contribution from the upper and lower sides $s=\sigma\pm iT$, $-c\leq\sigma\leq c'$, vanishes as $T\ra\infty$ provided $|\arg\,a|<\fs\pi$, since from the behaviour
\[\g(\sigma\pm it)=O(t^{\sigma-\fr}e^{-\fr\pi t}),\qquad \zeta(\sigma\pm it)=O(t^{\mu(\sigma)}\log^A t)\qquad(t\ra\infty),\]
where for $\sigma$ and $t$ real
\[\mu(\sigma)=0\ (\sigma>1),\quad \fs-\fs\sigma\ (0\leq\sigma\leq 1),\quad \fs-\sigma\ (\sigma<0),\]
\[A=1\ (0\leq\sigma\leq 1),\quad A=0\ \mbox{otherwise},\]
the modulus of the integrand is controlled by $O(T^{\sigma+\mu(\sigma)-\fr}\log\,T e^{-\Delta T})$, with $\Delta=\fs\pi-|\arg\,a|$. Displacement of the integration path to the right over a finite set of poles then yields (provided $w$ and $p$ are not even integers)
\bee\label{e22b}
S_p(a;w)=J_p(a;w)+\sum_{k=0}^{N-1}{}^{\!'}\frac{(-)^k}{k!}\,\zeta(w-kp)a^k+R_N,
\ee
where $N$ is a positive integer such that $N>s_0+\f{3}{2}$ and the prime on the sum over $k$ again denotes the omission of the term corresponding to $k=M$ when $w=pM+1$.

The remainder $R_N$ is given by
\begin{eqnarray*}
R_N&=&\frac{1}{2\pi i}\int_{c-\infty i}^{c+\infty i}\g(-s) \zeta(w-ps) a^{s}ds  \\
&=&\frac{(2\pi)^w}{2\pi i}\int_{c-\infty i}^{c+\infty i} \zeta(1-w+ps)\,\frac{\g(1-w+ps)}{\g(1+s)}\,\frac{\sin \fs\pi(ps-w)}{\sin \pi s}\,\frac{a^s}{(2\pi)^{ps}}ds
\end{eqnarray*}
where $c=N-\fs$ and the second expression follows from (\ref{e23a}).
Upon use of the result $|\zeta(x+iy)|<\zeta(x)$ when $x>1$, we obtain the bound
\bee\label{e23}
|R_N|<(2\pi)^{w-1} \zeta(pN\!-\!p(s_0\!-\!\fs)) \left(\frac{a}{(2\pi)^p}\right)^{\!N-\fr} \int_{-\infty}^\infty e^{-\phi t}F(t)\,dt,
\ee
where $\phi=\arg\,a$ and 
\[F(t)=\left|\frac{\g(1-w+ps)}{\g(1+s)}\right|\,\frac{\cosh \fs\pi pt}{\cosh \pi t}\qquad (s=N-\fs+it).\]
On the integration path, $F(t)$ is regular and satisfies $F(t)=O( e^{-\fr\pi |t|})$ as $t\ra\pm\infty$.
Hence the integral in (\ref{e23}) is convergent and independent of $|a|$ provided $|\phi|<\fs\pi$.
It then follows that 
\[\hspace{3cm}R_N=O(a^{N-\fr}) \qquad (a\ra 0 \ \mbox{in}\  |\arg\,a|<\fs\pi).\] 

The expansion (\ref{e22b}) is the dominant algebraic expansion associated with $S_p(a;w)$ valid as $a\ra 0$ in $|\arg\,a|<\fs\pi$, provided $w$ ($>0$) and $p$ are not even integers when the sum in (\ref{e22b}) is finite. The same analysis can be applied to the case with non-positive $w$ to yield the Berndt-Ramanujan result \cite[Theorem 3.1, p.~306]{BR} 
\[S_p(a;-w)=\frac{1}{p}\g\bl(\dfrac{1+w}{p}\br) a^{-(1+w)/p}+\sum_{k=0}^{N-1}\frac{(-)^k}{k!}\,\zeta(-w-kp)\,a^k+O(a^{N-\fr})\quad (w\geq 0)\]
as $a\ra 0$ in $|\arg\,a|<\fs\pi$. The reflection formula (\ref{e23a}) can be employed to convert the argument of the zeta function to a positive form.


\vspace{0.6cm}

\begin{center}
{\bf 3. \ The expansion of $S_p(a;w)$ when $w$ and $p$ are even integers}
\end{center}
\setcounter{section}{3}
\setcounter{equation}{0}
\renewcommand{\theequation}{\arabic{section}.\arabic{equation}}
Throughout this section we let $w$ and $p$ be even positive integers, with $w=2m$ where $m=1, 2, \ldots\,$. In this case $s_0=(2m-1)/p$, which cannot equal an integer and so no double pole can arise. More importantly, there is now only a finite set of poles of the integrand in (\ref{e21a})
at $s=s_0$ and $s=0, 1, 2, \ldots , K$, where $K=\lfloor w/p\rfloor$, since the poles of $\g(-s)$ at $s=K+k$ ($k=1, 2, \ldots$) are cancelled by the trivial zeros of the zeta function $\zeta(s)$ at $s=-2, -4, \ldots\,$. This has the consequence that the integrand is holomorphic in $\Re (s)>\max\{s_0,K\}$, so that further displacement of the contour can produce no additional algebraic terms in the expansion of $S_p(a;w)$. 

Thus, we find from (\ref{e22}) upon displacement of the integration path to the right over the poles of the integrand 
\bee\label{e31}
S_{p}(a;w)=\frac{1}{p}\g\bl(\frac{1-w}{p}\br) a^{(w-1)/p}+\sum_{k=0}^{K}\frac{(-)^k}{k!}\,\zeta(w-pk)\,a^k +(-)^m (2\pi)^w I_L,
\ee
where
\bee\label{e32}
I_L=\frac{(-)^m(2\pi)^{-w}}{2\pi i}\int_L \g(-s)\zeta(w-ps)a^{s}ds
\ee
and $L$ denotes a path parallel to the imaginary axis with $\Re (s)>(w/p)+\delta$, with $\delta$ denoting an arbitrary positive quantity.
This is easily seen to satisfy the requirement $\Re (s)>\max \{s_0,K_0\}$ necessary for the validity of (\ref{e31}). 
We now employ the functional relation for $\zeta(s)$ in (\ref{e23a}) 
to convert the argument of the zeta function in (\ref{e32}) into one with real part greater than unity. The integral in (\ref{e32}) can then be written in the form
\[I_L=\frac{1}{2\pi i}\int_L \zeta(1-w+ps)\,\frac{\g(1-w+ps)}{\g(1+s)}\,\frac{\sin \fs\pi ps}{\sin \pi s}\,\chi^{-s} ds,\]
where $\chi=(2\pi)^p a^{-1}$. 

In \cite[\S 8.1.4]{PK}, the zeta function appearing in the above integrand was written as an infinite series. Here we follow a suggestion made in \cite{B} and retain this function in the integrand; see also \cite[\S 8]{K}.
Making use of the expansion (see, for example, \cite[p.~368]{PK})
\[\frac{\sin \fs\pi ps}{\sin \pi s}=2\sum_{r=0}^{N-1} \cos \pi(\fs p\!-\!2r\!-\!1)s+\left\{\begin{array}{ll}\!\! 0 & (p/2\ \mbox{even})\\\!\! 1 & (p/2\ \mbox{odd})\end{array}\right.,\qquad N=[\f{1}{4} p],\]
where square brackets denote the {\it nearest\/} integer part\footnote{The nearest integer part corresponds to $[x]=N$ when $x$ is in the interval $(N-\frac{1}{2}, N+\frac{1}{2}]$. Note that when $p=2$, we have $N=0$ and the above expansion contains no information.}, we obtain 
\bee\label{e33a}
I_L=\sum_{r=0}^{N-1}\{J_{r}^++J_{r}^-\}+\left\{\begin{array}{ll}\!\! 0 & (p/2\ \mbox{even})\\\!\! J & (p/2\ \mbox{odd})\end{array}\right..
\ee
Here we have defined the integrals $J_r^\pm$ and $J$ by
\bee\label{e34}
J_{r}^\pm=\frac{1}{2\pi i}\int_L\zeta(1-w+ps)\,\frac{\g(1-w+ps)}{\g(1+s)}\,(\chi e^{\mp\pi i(p/2\!-\!2r\!-\!1)})^{-s}ds
\ee
and
\bee\label{e35}
J=\frac{1}{2\pi i}\int_L\zeta(1-w+ps)\,\frac{\g(1-w+ps)}{\g(1+s)}\,\chi^{-s}ds.
\ee
\vspace{0.3cm}

\noindent 3.1\ {\it Asymptotic evaluation of $J_r^\pm$ and $J$}
\vspace{0.2cm}

\noindent
The integrals $J_{r}^\pm$ and $J$ have no poles in the half-plane $\Re (s)>(w/p)+\delta$, so that we can displace the path $L$ as far to the right as we please. On such a displaced path $|s|$ is everywhere large. Let $M$ denote an arbitrary positive integer. The ratio of gamma functions appearing in (\ref{e34}) and (\ref{e35}) may then be expanded by making use of the result (for $p>1$) given in \cite[p.~53]{PK}
\bee\label{e36}
\frac{\g(1-w+ps)}{\g(1+s)}=\frac{A}{2\pi}(h\kappa^\kappa)^{-s}\bl\{\sum_{j=0}^{M-1} (-)^j c_j \g(\kappa s+\vartheta-j)+\rho_M(s) \g(\kappa s+\vartheta-M)\br\},
\ee
where $c_0=1$, $\rho_M(s)=O(1)$ as $|s|\ra\infty$ in $|\arg\,s|<\pi$ and
\bee\label{e36a}
\kappa=p-1,\qquad h=p^{-p},\qquad \vartheta=\fs-w,\qquad A=(2\pi)^\fr\,\kappa^{\fr-\vartheta} p^\vartheta.
\ee
The coefficients $c_j\equiv c_j(w,p)$ ($0\leq j\leq 4$) are listed in \cite[pp.~46--48]{PK} where an algorithm for their determination is described; see Section 4 for details.

Substitution of the expansion (\ref{e36}) into the integrals $J_{r}^\pm$ in (\ref{e34}) then produces
\begin{eqnarray}
J_{r}^\pm&=&\frac{A}{2\pi}\sum_{j=0}^{M-1} \frac{(-)^j c_j}{2\pi i}\int_L \zeta(1-w+ps)\,\g(\kappa s+\vartheta-j)\,(X e^{\mp\pi i\psi_r})^{-\kappa s} ds+{\cal R}_{M,r}^\pm \nonumber\\
&=&\frac{A}{2\pi\kappa}\sum_{j=0}^{M-1} (X e^{\mp\pi i\psi_r})^{\vartheta-j} 
\frac{(-)^j c_j}{2\pi i}\int_{L'} \g(u) \zeta(qu+\lambda_j) (Xe^{\mp\pi i\psi_r})^{-u}du+{\cal R}_{M,r}^\pm.\nonumber
\end{eqnarray}
Here we have made the change of variable $u\ra\kappa s+\vartheta-j$, with $L'$ denoting the modified integration path, and have defined
\bee\label{e313}
X:=\kappa(h\chi)^{1/\kappa},\quad \psi_r:=\frac{\fs p\!-\!2r\!-\!1}{\kappa}, \quad q:=\frac{p}{\kappa}, \quad \lambda_j:=1+\frac{1}{\kappa}(w+p(j-\fs))
\ee
together with the remainders
\bee\label{e37}
{\cal R}_{M,r}^\pm=\frac{A}{4\pi i}\int_L\rho_M(s)\zeta(1-w+ps)\,\g(\kappa s+\vartheta-M)\,(X e^{\mp\pi i\psi_r})^{-\kappa s} ds.
\ee
We note that $\lambda_j>0$ for $j\geq 0$ when $w>0$ and $p\geq 2$.
The above integrals  appearing in $J_r^\pm$ may now be evaluated by means of (\ref{e21a}), when we replace $s$ by $-s$
and allow the integration path $(c-\infty i, c+\infty i)$ to coincide with the path $L'$, to yield the sum $S_q(Xe^{\mp\pi i\psi_r};\lambda_j)$ as defined in (\ref{e11}). This evaluation is valid provided that
the variable $Xe^{\mp\pi i\psi_r}$ satisfies the convergence condition $|\arg (Xe^{\mp\pi i\psi_r})|<\fs\pi$; that is
\[\hspace{2cm}\left|\frac{\arg\,a}{\kappa}\mp\pi\psi_r\right|<\fs\pi \qquad (0\leq r\leq N-1).\]
It is routine to verify that these conditions are met when $|\arg\,a|<\fs\pi$. 

Thus we find
\bee\label{e37a}
J_r^\pm=\frac{A}{2\pi\kappa} \sum_{j=0}^{M-1} (-)^j c_j (Xe^{\mp\pi i\psi_r})^{\vartheta-j}\,S_q(Xe^{\mp\pi i\psi_r};\lambda_j)+{\cal R}_{M,r}^\pm.
\ee
Bounds for the remainders of the type ${\cal R}_{M,r}^\pm$ have been considered in \cite[p.~71, Lemma 2.7]{PK}; see also \cite[\S 10.1]{Br}. The integration path in (\ref{e37}) is such that $\Re (1-w+ps)>1$, so that we may employ the bound $|\zeta(x+iy)|\leq \zeta(x)$ for $x>1$. A slight modification of Lemma 2.7 in \cite[p.~71]{PK} then shows that
\bee\label{e37b}
{\cal R}_{M,r}^\pm=O\left(X^{\vartheta-M} e^{-Xe^{\mp\pi i\psi_r}}\right)
\ee
as $a\ra 0$ in the sector $|\arg\,a|<\fs\pi$.

An analogous procedure applied to $J$ in (\ref{e35}) shows that
\bee\label{e38}
J=\frac{A}{2\pi\kappa}\sum_{j=0}^{M-1} (-)^j c_jX^{\vartheta-j} S_q(X;\lambda_j)+O(X^{\vartheta-M}e^{-X})
\ee
as $a\ra 0$ in $|\arg\,a|<\fs\pi$.
\vspace{0.3cm}

\noindent 3.2\ {\it The expansion of $S_p(a;w)$}
\vspace{0.2cm}

\noindent
The expansion of $I_L$ as $a\ra 0$ in $|\arg\,a|<\fs\pi$ then follows from (\ref{e31}), (\ref{e33a}), (\ref{e37a}) and (\ref{e38}). We obtain the following theorem.
\newtheorem{theorem}{Theorem}
\begin{theorem}$\!\!\!.$
Let $m$ and $M$ be positive integers. Then, when $w=2m$ and $p$ is also an even positive integer, with $K=\lfloor w/p\rfloor$ and $N=[\f{1}{4}p]$ $($with square brackets denoting the nearest integer part$)$, we have the expansion valid as $a\ra 0$ in $|\arg\,a|<\fs\pi$ 
\bee\label{e310}
S_p(a;w)=\frac{1}{p}\g\bl(\frac{1-w}{p}\br) a^{(w-1)/p}+\sum_{k=0}^{K}\frac{(-)^k}{k!}\,\zeta(w-pk)\,a^k+ (-)^m(2\pi)^w I_L,
\ee
with 
\bee\label{e311}
I_L=\sum_{r=0}^{N-1} E_r(a;w,p)+\delta_{pp^*}\,{\hat E}_N(a;w,p),
\ee
where $\delta_{pp^*}$ is the Kronecker symbol with $p^*=4N+2$. 
The sums $E_r(a;w,p)$ are given by\footnote{The symbol $\sum_\pm$ signifies that the series with $\pm$ signs are to be added.} 
\bee\label{e312}
E_r(a;w,p)=\frac{A}{2\pi\kappa} \sum_{\pm}\sum_{j=0}^{M-1} (-)^jc_j (Xe^{\mp\pi i\psi_r})^{\vartheta-j} \,S_q(Xe^{\mp\pi i\psi_r};\lambda_j)+R_{M,r}
\ee
for $0\leq r\leq N-1$, where  
$X=\kappa (h(2\pi)^p/a)^{1/\kappa}$, $\psi_r=(\fs p\!-\!2r\!-\!1)/\kappa$ $(0\leq r\leq N-1)$, $\lambda_j=1+(w+p(j-\fs))/\kappa$, $q=p/\kappa$
and the parameters $\kappa$, $h$, $\vartheta$ and $A$ are defined in (\ref{e36a}). 
The leading coefficient $c_0=1$ and $c_j\equiv c_j(w,p)$ $(j\geq 1)$ are discussed in Section 4.    
The sum ${\hat E}_N(a;w,p)$ is also given by
(\ref{e310}) when we put $\psi_N\equiv 0$ and omit the summation $\sum_\pm$. The remainders $R_{M,r}$ satisfy the bound
\[R_{M,r}=O(\max\{X^{\vartheta-M} e^{-Xe^{\pm\pi i\psi_r}}\})\qquad (0\leq r\leq N).\]
\end{theorem}

It is seen from (\ref{e310}), (\ref{e311}) and (\ref{e312}) that the sum $S_p(a;w)$ has been expressed in terms of
the sums $S_q(Xe^{\mp\pi i\psi_r};\lambda_j)$, which involve the reciprocal power of the asymptotic variable $a$ scaling like $a^{-1/\kappa}$. Thus, as $a\ra 0$ the argument $X\ra\infty$. It is obvious from the definition in (\ref{e11}) (when $q>0$) that 
\[S_q(z;\lambda_j)\sim e^{-z}\qquad (z\ra\infty \ \mbox{in}\ |\arg\,z|<\fs\pi),\] 
so that the $E_r(a;w,p)$ represent a series of exponentially small expansions of increasing subdominance in the small-$a$ limit. In addition, the number of exponentially small expansions increases by one each time $p$ increases by 4. By means of a saddle-point analysis in the case $w=0$, this was demonstrated to correspond to a Stokes phenomenon when $p$ was allowed to vary continuously through the values $p=2, 6, 10, \ldots $; see \cite[\S\S 8.1.2, 8.1.7]{PK}. 
Finally, we remark that the exponents $p$ and $q$ are conjugate exponents \cite{B}, since 
\[\frac{1}{p}+\frac{1}{q}=1.\]

When $a$ is a real parameter, the expansion in Theorem 1 can be expressed in a different form by using (\ref{e11})
to represent the $S_q(Xe^{\mp\pi i\psi_r};\lambda_j)$ as infinite sums. Then from (\ref{e312})
we obtain the following theorem:
\begin{theorem}$\!\!\!.$
Let $w$ and $p$ be  even positive integers, $N=[\f{1}{4}p]$ and $M$ be a positive integer. Then, the exponentially small expansions in (\ref{e312}) valid as $a\ra 0+$ can be written in the form
\bee\label{e39b}
E_r(a;w,p)=\frac{A}{\pi\kappa} \sum_{n=1}^\infty n^{w-1}X_n^\vartheta\,e^{-X_n \cos \pi\psi_r}\,\Upsilon_{n,r} \qquad (0\leq r\leq N-1),
\ee
\bee\label{e39c}
{\hat E}_N(a;w,p)=\frac{A}{2\pi\kappa} \sum_{n=1}^\infty n^{w-1}X_n^\vartheta\,e^{-X_n}\,\Upsilon_{n,N}.
\ee 
The $\Upsilon_{n,r}$ $(0\leq r\leq N)$ have the asymptotic expansions
\bee\label{e39d}
\Upsilon_{n,r}=\sum_{j=0}^{M-1}(-)^jc_j X_n^{-j} \cos[X_n \sin \pi\psi_r+\pi(j-\vartheta) \psi_r]+O(X_n^{-M}),
\ee
where $X_n=Xn^{p/\kappa}=\kappa (h(2\pi n)^p/a)^{1/\kappa}$, $\psi_r=(\fs p-2r-1)/\kappa$, $\psi_N\equiv 0$
and the other quantities are as defined in Theorem 1.
\end{theorem}

The result in (\ref{e310}) and (\ref{e311}) is the analogue of the Poisson-Jacobi transformation in (\ref{e12}) corresponding to $w=0$, $p=2$. In this latter case, $N=0$ and, from (\ref{e311}), $I_L=E_0(a;0,2)$. The ratio of gamma functions in (\ref{e36}) is replaced by the single gamma function $\g(s+\fs)$ by the duplication formula for the gamma function, with the result that $c_0=1$, $c_j=0$ ($j\geq 1$) and consequently $\Upsilon_{n,N}=1$ for all $n\geq 1$. Then, since $K=0$ and $\zeta(0)=-\fs$, (\ref{e310}), (\ref{e311}) and (\ref{e39c}) reduce to (\ref{e12}). The resulting expansion is valid for all values of the parameter $a$ (not just $a\ra 0$) satisfying $|\arg\,a|<\fs\pi$. When $w=0$, $p=2m$, the expansions (\ref{e39b}) and (\ref{e39c}) reduce to those given in \cite{P}; see Section 5.

We observe that the $n$-dependence in the sums $E_r(a;w,p)$ from the factor
$n^{w-1} X_n^\vartheta$ is given by $n^{w-1}\,n^{p\vartheta/\kappa}=n^{-(2w+p-2)/(2\kappa)}$.
Since $p\geq 2$ and $w>0$, this is seen to correspond to a negative power of $n$.

\vspace{0.6cm}

\begin{center}
{\bf 4.\ The coefficients $c_j$}
\end{center}
\setcounter{section}{4}
\setcounter{equation}{0}
\renewcommand{\theequation}{\arabic{section}.\arabic{equation}}
We describe an algorithm for the computation of the coefficients $c_j\equiv c_j(w,p)$ that appear in the exponentially small expansions $E_r(a;w,p)$ in (\ref{e312}) and (\ref{e39d}).
The expression for the ratio of two gamma functions in (\ref{e36}), with $\alpha\equiv 1-w$ for convenience, takes the form
\[\frac{\g(\alpha+ps)}{\g(1+s) \g(\kappa s+\vartheta)}=\frac{A}{2\pi} (h\kappa^\kappa)^{-s}\bl\{\sum_{j=0}^{M-1} \frac{c_j}{(1-\kappa s-\vartheta)_j}+\frac{\rho_M(s)}{(1-\kappa s-\vartheta)_M}\br\},\]
where the parameters $\kappa$, $h$, $\vartheta$ and $A$ are defined in (\ref{e36a}) and $(\alpha)_j=\g(\alpha+j)/\g(\alpha)$ is the Pochhammer symbol. If we introduce the scaled gamma function $\g^*(z)=\g(z)/(\sqrt{2\pi}\,z^{z-\fr}e^{-z})$, 
then we have
\[\g(\beta s+\gamma)=\g^*(\beta s+\gamma) (2\pi)^\fr e^{-\beta s} (\beta s)^{\beta s+\gamma-\fr}\,{\bf e}(\beta s;\gamma),\]
where
\[{\bf e}(\beta s;\gamma):= \exp \bl[(\beta s+\gamma-\fs) \log (1+\frac{\gamma}{\beta s})-\gamma\br].\]
The above ratio of gamma functions may therefore be rewritten as
\bee\label{e41}
R(s)G(s)=\sum_{j=0}^{M-1}\frac{c_j}{(1-\kappa s-\vartheta)_j}+\frac{\rho_M(s)}{(1-\kappa s-\vartheta)_M}
\ee
as $|s|\ra\infty$ in $|\arg\,s|<\pi$, where
\[R(s)=\frac{{\bf e}(ps;\alpha)}{{\bf e}(s;1) {\bf e}(\kappa s;\vartheta)},\qquad G(s)=\frac{\g^*(\alpha+ps)}{\g^*(1+s) \g^*(\kappa s+\vartheta)}. \]

We now let $\xi:=(\kappa s)^{-1}$ and expand $R(s)$ and $G(s)$ for $\xi\ra 0$ making use of the well-known expansion \cite[p.~71]{PK}
\[\g^*(z)\sim\sum_{k=0}^\infty(-)^k\gamma_kz^{-k}\qquad(|z|\ra\infty;\ |\arg\,z|<\pi),\]
where $\gamma_k$ are the Stirling coefficients, with 
\[\gamma_0=1,\quad \gamma_1=-\f{1}{12},\quad \gamma_2=\f{1}{288},\quad  \gamma_3=\f{139}{51840},
\quad \gamma_4=-\f{571}{2488320}, \ldots\ .\]
After some straightforward algebra we find that
\[R(s)=1+\frac{\xi}{2}\bl\{\frac{\alpha(\alpha-1)\kappa}{p}-\vartheta(\vartheta-1)\br\}+O(\xi^2), \quad
G(s)=1+\frac{\xi}{12}\bl(1-p-\frac{1}{p}\br)+O(\xi^2),\]
so that upon equating coefficients of $\xi$ in (\ref{e41}) we can obtain $c_1$.
The higher coefficients can be obtained
by matching coefficients recursively with the aid of {\it Mathematica} to find \cite[p.~47]{PK}
\[c_0=1,\quad c_1=\frac{1}{24p}(2-5p+2p^2-12w+12pw+12w^2),\]
\[c_2=\frac{1}{1152p^2}(4+28p-87p^2+28p^3+4p^4+48w-216pw+24p^2w+144p^3w\]
\bee\label{e42}
\hspace{2cm}-96w^2-120pw^2+480p^2w^2-96w^3+480pw^3+144w^4), \ldots\ .
\ee
The rapidly increasing complexity of the coefficients with $j\geq 3$ prevents their presentation. However,
this procedure is found to work well in specific cases when the various parameters have numerical values, where up to a maximum of 100 coefficients have been so calculated.
In Table 1 we present the values\footnote{In the tables we write the values as $x(y)$ instead of $x\times 10^y$.} of the coefficients $c_j$ for $1\leq j\leq 8$ in the specific examples considered in Section 5.
\begin{table}[t]
\caption{\footnotesize{The coefficients $c_j$ ($1\leq j\leq 8$) for different $p$ and $w$.}}
\begin{center}
\begin{tabular}{|l|c|c||c|c|}
\hline
&&&&\\[-0.3cm]
\mcol{1}{|c|}{$j$} & \mcol{1}{c|}{$p=4,\ w=2$} &\mcol{1}{c||}{$p=4,\ w=4$}
& \mcol{1}{c|}{$p=6,\ w=2$} &\mcol{1}{c|}{$p=6,\ w=4$}\\
[.1cm]\hline
&&&&\\[-0.25cm]
1 & 1.395833\,$(0)$ & 3.645833\,$(0)$ & 1.472222\,$(0)$ & 3.305556\,$(0)$\\
2 & 3.495009\,$(0)$ & 1.648980\,$(1)$ & 3.861497\,$(0)$ & 1.469946\,$(1)$\\
3 & 1.230179\,$(1)$ & 9.075366\,$(1)$ & 1.380091\,$(1)$ & 8.081628\,$(1)$\\
4 & 5.555372\,$(1)$ & 5.899040\,$(2)$ & 6.207979\,$(1)$ & 5.260968\,$(2)$\\
5 & 3.060544\,$(2)$ & 4.424055\,$(3)$ & 3.387328\,$(2)$ & 3.949570\,$(3)$\\
6 & 1.990604\,$(3)$ & 3.760330\,$(4)$ & 2.188492\,$(3)$ & 3.358058\,$(4)$\\
7 & 1.493190\,$(4)$ & 3.572267\,$(5)$ & 1.639364\,$(4)$ & 3.189927\,$(5)$\\
8 & 1.269216\,$(5)$ & 3.750863\,$(6)$ & 1.396172\,$(5)$ & 3.348999\,$(6)$\\
[.2cm]\hline
\end{tabular}
\end{center}
\end{table}

When $p=2$, use of the duplication formula shows that the ratio of gamma functions in (\ref{e36})
becomes 
\[\frac{\g(\fs\!-\!\fs w\!+\!s) \g(1\!-\!\fs w\!+\!s)}{\g(1+s)}=\sum_{j=0}^{M-1} (-)^j c_j \g(s\!+\!\fs\!-\!w\!-\!j)+\rho_M(s) \g(s\!+\!\fs\!-\!w\!-\!M).\]
The coefficients $c_j$ in this case can be expressed in closed form as \cite[p.~53]{PK}
\bee\label{e43}
c_j=\frac{(\fs w)_j (\fs+\fs w)_j}{j!}=\frac{2^{-2j}(w)_{2j}}{j!} \qquad (p=2).
\ee
Finally, we mention that when $w=0$ (corresponding to the Euler-Jacobi series) the coefficients $c_j$ are listed  for $p\geq 2$ and $j\leq 8$ in \cite[p.~374]{PK}.
\vspace{0.6cm}

\begin{center}
{\bf 5.\ Numerical results and concluding remarks}
\end{center}
\setcounter{section}{5}
\setcounter{equation}{0}
\renewcommand{\theequation}{\arabic{section}.\arabic{equation}}
We present some examples of the expansion of $S_p(a;w)$ given in Theorem 1 when $p$ and $w=2m$ are even integers.
For convenience in presentation, we extract the factor $e^{-z}$ from the sum $S_q(z;\lambda_j)$ by writing
\[S_q(z;\lambda_j)=e^{-z} {\hat S}_2(z;\lambda_j),\qquad {\hat S}_q(z;\lambda_j):=\sum_{n=1}^\infty \frac{e^{-z(n^q-1)}}{n^{\lambda_j}}.\]
It follows that, when $q>0$, ${\hat S}_q(z;\lambda_j)=O(1)$ as $z\ra\infty$ in $|\arg\,z|<\fs\pi$.
\vspace{0.3cm}

\noindent{\it Example 1.}\ \ \ In the case $p=2$, we have $\kappa=1$, $q=2$, $N=0$, $K=m$, $p^*=2$, $\psi_0\equiv 0$ and $X=\pi^2/a$. The quantity $\delta_{pp^*}=1$ so that the exponentially small component of $S_2(a;2m)$ consists of the single term ${\hat E}_0(a;2m,2)$. From  (\ref{e310}), (\ref{e311}), (\ref{e312}) and (\ref{e43}) we therefore find 
\[S_2(a;2m)-\frac{1}{2}\g\bl(\frac{1-2m}{2}\br) a^{m-\fr}-\sum_{k=0}^m\frac{(-)^k}{k!}\,\zeta(2m-2k) a^k\]
\bee\label{e51}
=(-)^m \left(\frac{a}{\pi}\right)^{\!2m-\fr}e^{-\pi^2/a} \bl\{\sum_{j=0}^{M-1}\frac{(-)^j(2m)_{2j}}{j!} \left(\frac{a}{4\pi^2}\right)^j \!{\hat S}_2(\pi^2/a;2m+2j)+O(a^M) \br\}
\ee
as $a\ra 0$ in the sector $|\arg\,a|<\fs\pi$. The expansion in this case has been given in an equivalent form in \cite{P}.
\vspace{0.3cm}

\noindent{\it Example 2.}\ \ \ 
When $p=4$, we have $\kappa=3$, $q=\f{4}{3}$, $N=1$, $K=\lfloor\fs m\rfloor$, $p^*=6$ and
\[X=3(\fs\pi)^{4/3}\,a^{-1/3},\qquad \lambda_j=\f{1}{3}(2m+4j+1).\]
The quantity $\delta_{pp^*}=0$ so that there is the single exponentially small term $E_0(a;2m,4)$ with $\psi_0=\f{1}{3}$. Then we find the expansion as $a\ra 0$ in $|\arg\,a|<\fs\pi$ given by
\[S_4(a;2m)
-\frac{1}{4}\g\bl(\frac{1-2m}{4}\br) a^{(2m-1)/4}-\sum_{k=0}^K\frac{(-)^k}{k!}\,\zeta(2m-4k)a^k\]
\bee\label{e52}
=(-)^m\bl(\frac{2a}{\pi}\br)^{\!(4m-1)/6}\sum_\pm e^{-Xe^{\mp\frac{1}{3}\pi i}\mp\frac{1}{3}\pi i\vartheta}
\bl\{\sum_{j=0}^{M-1}\frac{(-)^j c_j}{(Xe^{\mp\frac{1}{3}\pi i})^j}\,{\hat S}_\frac{4}{3}(Xe^{\mp\frac{1}{3}\pi i};\lambda_j)+O(a^{M/3})\br\}
\ee
where the coefficients $c_j\equiv c_j(2m,4)$ can be obtained from (\ref{e42}) and \cite[p.~47]{PK} as
\[c_0=1,\quad c_1=\frac{1}{48} (7+36m+24m^2),\]
\[c_2=\frac{1}{4608} (385+4392m+7104m^2+3648m^3+576m^4),\]
\[c_3=\frac{1}{663552}(39655+1191132m+2970936m^2+2666880m^3\]\[\hspace{5cm}+1080000m^4+200448m^5+13824m^6), \ldots \ .\]
These coefficients are listed in Table 1 for $1\leq j\leq 8$ when $m=1$ and $m=2$.
\vspace{0.3cm}

\noindent{\it Example 3}.\ \ 
When $p=6$, we have $\kappa=5$, $q=\f{6}{5}$, $N=1$, $K=\lfloor \f{1}{3}m\rfloor$ and 
\[X=5(\f{1}{3}\pi)^{6/5}\,a^{-1/5},\qquad \lambda_j=\f{1}{5}(2m+6j+2).\]
In this case $p^*=6$, so that $\delta_{pp^*}=1$ and there are now two exponentially small expansions $E_0(a;2m,6)$, with $\psi_0=\f{2}{5}$, and ${\hat E}_1(a;2m,6)$. Then, as $a\ra 0$ in $|\arg\,a|<\fs\pi$, we have the expansion
\[S_6(a;2m)-\frac{1}{6}\g\bl(\frac{1-2m}{6}\br) a^{(2m-1)/6}-\sum_{k=0}^K\frac{(-)^k}{k!}\,\zeta(2m-6k)a^k\]
\[=(-)^m\bl(\frac{3a}{\pi}\br)^{\!(4m-1)/10}\,\sum_\pm e^{-Xe^{\mp\frac{2}{5}\pi i}\mp\frac{2}{5}\pi i\vartheta}
\bl\{\sum_{j=0}^{M-1}\frac{(-)^j c_j}{(Xe^{\mp\frac{2}{5}\pi i})^j}\,{\hat S}_\frac{6}{5}(Xe^{\mp\frac{2}{5}\pi i};\lambda_j)+O(a^{M/5})\br\}\]
\bee\label{e53}
+(-)^m\bl(\frac{3a}{\pi}\br)^{\!(4m-1)/10}\,e^{-X}
\bl\{\sum_{j=0}^{M-1}\frac{(-)^j c_j}{X^j}\,{\hat S}_\frac{6}{5}(X;\lambda_j)+O(a^{M/5})\br\}.
\ee
The first few coefficients $c_j\equiv c_j(2m,6)$ are
\[c_0=1,\quad c_1=\frac{1}{36}(11+30m+12m^2),\]
\[c_2=\frac{1}{2592} (517+3840m+4116m^2+1392m^3+144m^4),\]
\[c_3=\frac{1}{1399680}(-22253+426550m+8181720m^2+5237640m^3\]\[\hspace{5cm}+1468800m^4+185760m^5+8640m^6).\]
These coefficients are listed in Table 1 for $1\leq j\leq 8$ when $m=1$ and $m=2$.
\vspace{0.3cm}


\begin{table}[th]
\caption{\footnotesize{Values of the absolute error in the computation of $S_p(a;w)$ defined in (\ref{e54}) using the expansions (\ref{e51}) and (\ref{e52}).
The value of the index $j_0$ corresponds to optimal truncation of the subdominant expansion $E_0(a;w,p)$.}}
\begin{center}
\begin{tabular}{|l|llr|llr|}
\hline
&&&&&&\\[-0.25cm]
\mcol{1}{|c|}{} & \mcol{3}{c|}{$p=2,\ w=2$} & \mcol{3}{c|}{$p=2,\ w=4$} \\
\mcol{1}{|c|}{$a$} & \mcol{1}{c}{$|{\cal S}_{p,w}|$} & \mcol{1}{c}{$|{\cal S}_{p,w}-E_0|$}  & \mcol{1}{c|}{$j_0$} & \mcol{1}{c}{$|{\cal S}_{p,w}|$} & \mcol{1}{c}{$|{\cal S}_{p,w}-E_0|$} & \mcol{1}{c|}{$j_0$}\\
[.1cm]\hline
&&&&&&\\[-0.3cm]
1.00 & $8.146(-06)$ & $6.637(-09)$ & 8  & $6.252(-07)$ & $3.642(-08)$ &  6 \\
0.75 & $2.031(-07)$ & $8.089(-12)$ & 11 & $9.296(-09)$ & $4.659(-11)$ &  9 \\
0.50 & $1.584(-10)$ & $1.260(-17)$ & 18 & $3.437(-12)$ & $7.635(-17)$ & 16 \\
0.20 & $5.774(-24)$ & $1.542(-43)$ & 47 & $2.189(-26)$ & $9.830(-43)$ & 45 \\
0.10 & $7.667(-46)$ & $1.486(-86)$ & 97 & $7.506(-49)$ & $9.631(-86)$ & 95 \\ 
[.2cm]\hline
&&&&&&\\[-0.25cm]
\mcol{1}{|c|}{} & \mcol{3}{c|}{$p=4,\ w=2$} & \mcol{3}{c|}{$p=4,\ w=4$} \\
\mcol{1}{|c|}{$a$} & \mcol{1}{c}{$|{\cal S}_{p,w}|$} & \mcol{1}{c}{$|{\cal S}_{p,w}-E_0|$}  & \mcol{1}{c|}{$j_0$} & \mcol{1}{c}{$|{\cal S}_{p,w}|$} & \mcol{1}{c}{$|{\cal S}_{p,w}-E_0|$} & \mcol{1}{c|}{$j_0$}\\
[.1cm]\hline
&&&&&&\\[-0.3cm]
0.200 & $3.473(-03)$ & $1.329(-06)$ & 7  & $3.919(-04)$ & $8.742(-06)$ &  6 \\
0.100 & $4.863(-04)$ & $2.749(-08)$ & 11 & $4.805(-05)$ & $4.879(-07)$ &  8 \\
0.050 & $2.737(-05)$ & $2.156(-10)$ & 14 & $8.456(-06)$ & $1.420(-09)$ & 11 \\
0.010 & $4.221(-09)$ & $4.621(-17)$ & 23 & $7.982(-09)$ & $3.041(-16)$ & 21 \\
0.001 & $1.064(-14)$ & $1.033(-36)$ & 53 & $1.876(-16)$ & $6.799(-36)$ & 51 \\ 
[.2cm]\hline
\end{tabular}
\end{center}
\end{table}

We show the results of numerical calculations to demonstrate the achievable accuracy of the expansion in Theorem 1.
We define the difference between $S_p(a;w)$ and the finite algebraic expansion by
\bee\label{e54}
{\cal S}_{p,w}\equiv {\cal S}_{p,w}(a):=S_p(a;w)-\frac{1}{p}\g\bl(\frac{1-w}{p}\br) a^{(w-1)/p}-\sum_{k=0}^{K}\frac{(-)^k}{k!}\,\zeta(w-pk)\,a^k.\ee
In Table 2 we present the absolute error in the computation of $S_p(a;w)$ for different values of the parameter $a$ in the two cases $p=2$ and $p=4$, with $w=2$ and $w=4$ using the expansions given in (\ref{e51}) and (\ref{e52}).
The first column in each entry displays the absolute value of ${\cal S}_{p,w}$; that is, the accuracy achievable with just the algebraic expansion and no subdominant exponential terms. The second column shows the absolute error when the single {\it optimally truncated\/}  exponential expansion $E_0(a;w,p)$ (denoted by $E_0$ in the table) is included. The optimal truncation index $j_0$, corresponding to truncation of the exponential expansion $E_0(a;w,p)$ at, or near, the least term in magnitude, is indicated in the final column.

The situation when there is only a single subdominant exponentially small expansion present is straightforward:
this sum is truncated at some suitable point thereby introducing a truncation error. If truncation is optimal, then the resulting error is exponentially more recessive than the parent exponential expansion. However, in the case of two, or more, exponential expansions of different degrees of subdominance (corresponding to $p\geq 6$) the situation is not so obvious. It is not clear, without further investigation, how the error from the truncated leading exponential series compares with the contribution from the next series. 

\begin{table}[h]
\caption{\footnotesize{Values of the absolute error in the computation of $S_6(a;2)$ defined by (\ref{e54}) using the expansion (\ref{e53}). The value of the index $j_0$ corresponds to optimal truncation of the expansion $E_0(a;2,6)$.}}
\begin{center}
\begin{tabular}{|l|llrc|ll|}
\hline
&&&&&&\\[-0.25cm]
\mcol{1}{|c|}{$a$} & \mcol{1}{c}{$|{\cal S}_{6,2}|$} & \mcol{1}{c}{$|{\cal S}_{6,2}-E_0|$}  & \mcol{1}{c}{$j_0$} & \mcol{1}{c|}{$|{\cal S}_{6,2}-E_{0,1}|$} & \mcol{1}{c}{Min $|E_0|$} & \mcol{1}{c|}{$E_1(j=0)$}\\
[.1cm]\hline
&&&&&&\\[-0.3cm]
1$\times 10^{-1}$ & $2.935(-02)$ & $3.780(-05)$ & 6  & $-\!\!-$     & $9.422(-05)$ & $5.095(-05)$ \\
5$\times 10^{-2}$ & $1.617(-03)$ & $3.037(-05)$ & 8  & $1.200(-05)$ & $1.729(-05)$ & $1.191(-05)$ \\
1$\times 10^{-2}$ & $9.512(-04)$ & $1.193(-07)$ & 12 & $5.339(-08)$ & $1.228(-07)$ & $1.904(-07)$ \\
5$\times 10^{-3}$ & $1.292(-03)$ & $1.099(-08)$ & 13 & $8.713(-09)$ & $9.090(-09)$ & $2.148(-08)$ \\
1$\times 10^{-3}$ & $1.604(-04)$ & $3.452(-11)$ & 19 & $3.483(-12)$ & $3.757(-12)$ & $4.053(-11)$ \\
1$\times 10^{-4}$ & $9.894(-07)$ & $8.801(-17)$ & 31 & $2.230(-19)$ & $3.024(-19)$ & $9.201(-17)$ \\ 
1$\times 10^{-5}$ & $6.209(-10)$ & $1.522(-25)$ & 51 & $1.963(-30)$ & $1.964(-30)$ & $1.564(-25)$ \\
[.2cm]\hline
\end{tabular}
\end{center}
\end{table}

We illustrate this by considering the case $p=6$ and $w=2$ given in (\ref{e53}). In Table 3 we present the absolute error in the computation of $S_6(a;2)$ as a function of the parameter $a$. We show, in order, the value of $|{\cal S}_{6,2}|$ and the absolute error in ${\cal S}_{6,2}-E_0(a;2,6)$ when the leading subdominant exponential expansion $E_0(a;2,6)$ is optimally truncated at index $j_0$. The fourth column gives the absolute error when the first few terms of the second exponential expansion ${\hat E}_1(a;2,6)$ are included (for brevity in the table these exponential expansions are labelled $E_0$ and $E_1$, and their sum is denoted by $E_{0,1}$). The final two columns show the values of the least term (including prefactors) in
$E_0(a;2,6)$ at optimal truncation and the values of the leading term ($j=0$) of the sub-subdominant expansion ${\hat E}_1(a;2,6)$.

A cursory inspection of Table 3 shows that for $a\simeq 0.1$ the leading term of ${\hat E}_1(a;2,6)$ is less than the minimum term of $E_0(a;2,6)$ and consequently that inclusion of ${\hat E}_1(a;2,6)$ cannot improve the accuracy.
For $a\,\ltwid\,0.01$, the reverse is true: the leading terms of ${\hat E}_1(a;2,6)$ are greater than the minimum term of $E_0(a;2,6)$ and their inclusion therefore increases the overall accuracy. However, it is clear that in both cases the final accuracy achievable is limited by the optimal truncation of the leading subdominant expansion $E_0(a;2,6)$. Further improvement in the accuracy would require a hyperasymptotic treatment in order to deal with the divergent tails of $E_0(a;2,6)$ and ${\hat E}_1(a;2,6)$. A possible hyperasymptotic scheme for the Euler-Jacobi series with $p=3$ and $w=0$ has been discussed in \cite[\S 8]{K}.

Finally we remark that the asymptotics of the alternating version of (\ref{e11}) can be deduced from the result in Theorem 1 by making use of the identity
\[\sum_{n=1}^\infty (-)^n\frac{e^{-an^p}}{n^w}=2^{1-w} S_p(2^pa;w)-S_p(a;w).\]

\end{document}